\documentclass[12pt]{amsart}
\usepackage{amscd,amssymb}
\theoremstyle{plain}
\newtheorem{thm}[subsection]{Theorem}
\newtheorem{lem}[subsection]{Lemma}
\newtheorem{prop}[subsection]{Proposition}

\theoremstyle{definition}
\newtheorem{rem}[subsection]{Remark}
\newtheorem{defn}[subsection]{Definition}
\newtheorem{exm}[subsection]{Example}

\numberwithin{equation}{section}

\newcommand{\Z}{{\mathbb Z}}
\newcommand{\Q}{{\mathbb Q}}

\newcommand{\C}{{\mathbb C}}
\renewcommand{\H}{{\mathbb H}}
\newcommand{\Ca}{{\mathbf{Ca}}}
\newcommand{\F}{{\mathbb F}}
\renewcommand{\P}{{\mathbb P}}

\newcommand{\cf}{{\it cf.\ \/}}
\newcommand{\eg}{{\it e.g.\ \/}}

\renewcommand{\a}{{\alpha }}

\newcommand{\mero}{{ \, - -\negmedspace \to\, }}

\DeclareMathOperator{\vol}{{vol}}

\DeclareMathOperator{\sys}{{sys}}

\DeclareMathOperator{\tors}{{tors}}
\DeclareMathOperator{\rank}{{rank}}
\DeclareMathOperator{\id}{{id}}

\begin{document}

\title[Volumes, systoles, and Whitehead products] %
{Volumes, middle-dimensional systoles, and Whitehead products}

\author{Ivan K.~Babenko}
\address{D\'epartement des Sciences Math\'ematiques,
Universit\'e de Montpellier 2, place Eug\`ene Bataillon,
34095 Montpellier, France}
\email{babenko@darboux.math.univ-montp2.fr}

\author{Mikhail G.~Katz}
\address{UMR 9973, D\'epartement de Math\'ematiques,
Universit\'e de Nancy 1, B.P. 239, 54506 Vandoeuvre, France}
\email{katz@iecn.u-nancy.fr}

\author{Alexander I.~Suciu}
\address{Department of Mathematics,
Northeastern University,
Boston, MA 02115, USA}
\email{alexsuciu@neu.edu}

\thanks{The second author is grateful to the Geometry \& Dynamics
Program at Tel Aviv University for hospitality during part of
the preparation of this work.}
\thanks{The third author was partially supported by
N.S.F.~grant DMS--9504833.  He is grateful to the
Institut \'Elie~Cartan at Henri~Poincar\'e University--Nancy~1
for hospitality during the completion of this work.}


\begin{abstract}  
Let $X$ be a closed manifold of dimension $2m\ge 6$ with torsion-free
middle-dimensional homology.  We construct metrics on $X$ of
arbitrarily small volume, such that every middle-dimensional
submanifold of less than unit volume necessarily bounds.  Thus,
Loewner's theorem has no higher-dimensional analogue.
\end{abstract}

\maketitle

\section{Introduction}
\label{sec:intro}

Let $(X,g)$ be a closed, orientable Riemannian manifold of even
dimension $2m$.  The following notion was introduced by Marcel~Berger
in~\cite{Ber1}, \cite{Ber2}.

\begin{defn} \label{def:systole}
The {\it $k$-systole} of $(X,g)$, $\sys_k(g)$, is the infimum of areas
of non-bounding cycles represented by maps of $k$-dimensional
manifolds into $X$.
\end{defn}

In this note, we will be interested in the following question: Does
there exist a constant, $C$, such that every metric $g$ on $X$
satisfies
\begin{equation} \label{eq:question}
\sys_{m}^{2}(g) \le C\cdot \vol_{2m}(g) ?
\end{equation}
If there is no such $C$, we say that $X$ is {\em systolically free}.

In the case of surfaces of positive genus the answer to question
\eqref{eq:question} is affirmative.  For a history of the problem
(dating to C. Loewner's solution in the case $X=S^1\times S^1$), see
M.~Berger \cite{Ber2} and P.~Sarnak~\cite{Sa}.  In the case $m\ge 2$,
this question has been referred to by M. Gromov as the ``basic
systolic problem" (\cite{G2}, p.~357); see also the subsection
``Systolic reminiscences" of \cite{G3}, p.~267.  Progress on the
problem became possible once Gromov described a special family of
metrics on $S^1\times S^3$, and surgical procedures suitable for
generalizations (\cite{G3}, p.~264); see also \cite{K1} and \cite{BaK}.

The purpose of this note is to prove the following result
(see Theorem~\ref{thm:mainCW} for a statement in the more
general context of CW-complexes).

\begin{thm} \label{thm:main}
Let $X$ be a closed, orientable, smooth manifold of dimension $2m$.
If $m\ge 3$ and $H_m(X)$ is torsion-free, then $X$ is systolically free.
\end{thm}

An underlying theme of this paper is the influence of homotopy theory
on the geometric inequality \eqref{eq:question}.  The basic idea is to
establish the systolic freedom of a complicated manifold, $X$, by
``folding away" some of the topology of $X$ into lower-dimensional
cells attached to a simpler manifold, $Y$, whose systolic freedom is
already established.  The notion of ``folding away" is captured in
Definition~\ref{def:mero} of a ``meromorphic map," $f: X \mero Y$.
The starting point is Gromov's sequence of metrics \eqref{eq:gromov},
which establishes the systolic freedom of products of spheres.

Our topological tools are, on one hand, the Hilton-Milnor theorem
calculating homotopy groups of a wedge of spheres, and, on the other hand,
theorems of B.~Eckmann and G.~Whitehead on composition maps in homotopy
groups of spheres.  Our geometric tools are the coarea inequality
of \cite{G1} and pullback arguments for simplicial metrics as described
by the first author in~\cite{Bk}.

The case $m=2$ of \eqref{eq:question} remains open, but it has been reduced
to either $\C\P^2$ or $S^2\times S^2$ by the first two authors in~\cite{BaK}.
Even if we restrict the class of competing metrics to homogeneous
ones, the inequality \eqref{eq:question} is violated in certain
cases such as $S^3 \times S^3$, see~\cite{K2}.

One may view our construction as a way of producing metrics for which
the systole and the mass in middle dimension do not agree,
compare~\cite{Gluck}.

An announcement of this paper appeared in~\cite{K2}.

The structure of the paper is as follows:

\begin{itemize}
\item
In section~\ref{sec:mero}, we define ``meromorphic maps" between regular
CW-complexes, which allows us to correlate their systolic freedom.
We use this technique to give a short proof of the systolic freedom
of $S^m \times S^m$.
\item
In section~\ref{sec:homotopy},
we find maps from the $(2m-1)$-skeleton of $X$
to a wedge of $m$-spheres that induce monomorphisms in
$H_m(-,\Q)$, and self-maps of $\vee S^m$ that send
$\pi_{2m-1} (\vee S^m)$ to the subgroup generated
by Whitehead products.
\item In section~\ref{sec:b1}, we prove our theorem in the case where
  $b_m(X)=1$, by mapping $X$ meromorphically to $S^m \times S^m$.
\item In section~\ref{sec:bgt1}, we present the proof in the general
  case.  This is achieved by mapping $X$ meromorphically to the
  $2m$-skeleton of a product of sufficiently many $m$-spheres.
\end{itemize}

\section{Systolic freedom of CW-complexes and meromorphic maps}
\label{sec:mero}

In order to prove Theorem~\ref{thm:main}, we will enlarge the
class of manifolds to that of piecewise smooth, simplicial complexes,
for which one can still define metrics, volumes, and systoles.
We will actually prove our theorem in the context of
finite, regular CW-complexes.
Such a complex $K$ can be triangulated so that the resulting
simplicial complex is a subdivision of $K$ (see \cite{LW}, p.~80).

\begin{defn} \label{def:free}
A finite, regular CW-complex $K$ of dimension $2m$
is called {\it systolically free} if
\begin{equation}  \label{eq:free}
\inf_{g} \frac{\vol_{2m}(g)}{\sys_m^2(g)}=0,
\end{equation}
where the infimum is taken over all metrics $g$ on $K$.
This amounts to the existence of a sequence of metrics $\{g_j\}$
such that \begin{equation}  \label{eq:freeseq}
\sys_m^2(g_j)\ge j\, \vol_{2m}(g_j).
\end{equation}
\end{defn}

\begin{rem}  \label{rem:regCW}
The systolic freedom of $K$ (or the absence thereof) is
independent of the piecewise smooth simplicial structure
that one chooses in its homotopy type.  This independence is
verified by means of the simplicial approximation theorem and
by the pullback arguments for metrics from \cite{Bk}, \cite{BaK}.
\end{rem}

\begin{thm} \label{thm:mainCW}
Let $K$ be a finite, regular CW-complex of dimension $2m\ge 6$.
If $H_m(K)$ is torsion-free, then $K$ is systolically free.
\end{thm}

This theorem (which slightly generalizes Theorem~\ref{thm:main}),
will be proved at the end of section~\ref{sec:bgt1}.
The key to the proof is the following notion,
inspired by complex analysis and surgery theory.

\begin{defn}  \label{def:mero}
Let $X$ and $Y$ be $2m$-dimensional CW-complexes.
A {\it ``meromorphic map''} from $X$ to $Y$ is a continuous map
$f:X\to W$ such that
\begin{enumerate}
\item $W$ is a CW-complex obtained from $Y$ by
attaching cells of dimension at most $2m-1$; \label{mero1}
\item $f_*: H_m(X) \to H_m(W)$ is a monomorphism.  \label{mero2}
\end{enumerate}
\end{defn}

We shall denote such ``maps'' by $f: X \mero Y$, and
drop the quotation marks.

\begin{exm}  \label{exm:blowup}
Let $X$ be a complex surface and $\widehat X \rightarrow X$
its blow-up at a point $p\in X$.  Then the classical
meromorphic map $X\rightarrow \widehat{X}$ can be modified in a
neighborhood of $p$ and extended to a continuous map from $X$ to
$\widehat{X}\cup_f\! B^3$ where the $3$-ball is attached along the
exceptional curve.
\end{exm}

\begin{prop} [\cite{BaK}]  \label{pull}
Let $X$ and $Y$ be finite, regular CW-complexes.
Suppose $X$ admits a meromorphic map to $Y$.  If\, $Y$
is systolically free, then $X$ is also systolically free.
\end{prop}

\begin{proof} [Proof (sketch)]
Let $f: X \to W=Y\cup \bigcup_i B^{k_i}_{i}$ be the given meromorphic map.
The attached cells (of dimension $k_i<2m$) do not affect the
$2m$-dimensional volume, and thus $W$ is still systolically free,
by the cylinder construction of \cite{BaK}, Lemma~6.1.
We now pull back to $X$ the systolically free metrics on $W$.
Thus $X$ is systolically free.
\end{proof}

\begin{prop} \label{exm:surgery}
  Let $X$ be an orientable, smooth manifold of dimension $n=2m$.
  Suppose $Y$ is obtained from $X$ by performing surgery on embedded,
  framed $k$-spheres, with $1\le k<m$.  Then there exists a
  meromorphic map from $X$ to $Y$.
\end{prop}
\begin{proof}
Let $W^{n+1}$ be the cobordism between $X$ and $Y$ defined
by the surgeries.  We claim that the inclusion $X\hookrightarrow W$
is the desired meromorphic map. Indeed, $W$ is obtained by attaching
handles $D^{k+1}\times D^{n-k}$ to $X\times I$, or dually,
by attaching handles $D^{n-k}\times D^{k+1}$ to $Y\times I$.
Hence, $W$ has the homotopy type of $Y$, with cells of
dimension $n-k\le n-1$ attached to it, and so
condition~(\ref{mero1}) is satisfied.   Since $k<m$,
condition~(\ref{mero2}) is satisfied also.
\end{proof}

The main geometric ingredient in the proof of Theorem \ref{thm:mainCW}
is the special case of a product of spheres, first proved in \cite{K1}.
We provide a different proof in Proposition~\ref{prop:snsn} below,
using CW-complexes and meromorphic maps.

In order to give the reader a heuristic understanding of the fundamental
examples, let us begin by describing Gromov's construction of systolically
free metrics on $S^1\times S^3$, see~\cite{G3}, paragraph~4.45.  We translate
Gromov's succinct example, from the language of global Riemannian
geometry into the dual language of differential forms.

Let $S^1$ be the unit circle with standard $1$-form $dz$,
and let $S^3$ be the unit $3$-sphere with standard contact $1$-form $b$.
Gromov's metrics $\{g_j\}$ are obtained by modifying the product metric
of a circle of length $2\pi$ with a $3$-sphere of radius $R=\sqrt{1+j^2}$,
by adding a non-diagonal term $-2jb\, dz$ (symmetric tensor product).
Explicitly, let us complete $b$ to a basis $(b,b',b'')$ of $1$-forms
which is orthonormal with respect to the metric of unit radius, so that
$db=b'\wedge b''$.  Then Gromov's metrics are given by:
\begin{equation}  \label{eq:gromov}
g_j=(dz-jb)^2+{b}^2+R^2({b'}^2+{b''}^2).
\end{equation}

Roughly speaking, from the point of view of an individual Hopf fiber
of $S^3$, each metric looks like a square torus of size $2\pi$ with
respect to a judiciously chosen basis.  Meanwhile, from the point
of view of the hypersurface $S^3$ of radius $R$, the metric $g_j$
looks like a family of metrics on $S^3$ parametrized by a circle of
length $\frac{2\pi}{R}$.

Let $*$ be the Hodge star operator of the metric $g_j$.  Then:
\begin{equation} \label{eq:star}
\begin{split}
*(dz)=*(dz-jb)+j*b= (b-j(dz-jb)) \wedge R^2 db \\
\qquad = R^4 b\wedge db - jR^2 dz\wedge db.
\end{split}
\end{equation}
The form $*dz$ is closed, since $d*dz=R^4 db\wedge db=0$.  We
normalize $*dz$ to obtain a calibrating form $\frac{1}{R}*dz$, whose
restriction to the $3$-sphere $S^3\subset S^1\times S^3$ coincides with
the volume form $R^3b\wedge db$ of this round sphere of radius $R$.
Estimating the volume and the systoles yields:
\begin{equation}  \label{lim1}
\lim_{j\to\infty} \frac{\vol_{4} (g_j)}{\sys_1(g_j)\cdot \sys_3(g_j)} =
\lim_{j\to\infty} \frac{j^2}{1\cdot j^3} = 0.
\end{equation}

This example generalizes to all products of spheres, except
$S^1\times S^1$ and possibly $S^2\times S^2$,
see \cite{K1}, \cite{BeK}, \cite{P}, \cite{BaK}.
For each such product $S^{m}\times S^{k}$, there is a sequence
of metrics $\{g_{j}\}$ such that
\begin{equation}  \label{systseq}
\lim_{j\to\infty} \frac{\vol_{m+k} (g_j)}{\sys_k(g_j)\sys_m(g_j)} =0.
\end{equation}
The existence of metrics on $S^m\times S^k$ satisfying \eqref{systseq}
is proved in two steps, first for $m>k$, and next for $m=k$.
The first step is done by direct geometric construction,
see \cite{BaK}, proof of Proposition~4.2.  The second step
was actually done beforehand, in \cite{K1}.  In order
to give a flavor of the arguments involved, we provide a
self-contained proof of the second step, assuming the first step.

\begin{prop}[\cite{K1}] \label{prop:snsn}
For $m\ge 3$, the product $S^m\times S^{m}$ is systolically free.
\end{prop}
\begin{proof}
Choose an integer $k$ such that $1\le k <m$.
Performing surgery on a standard $S^k\subset S^m$ yields
$S^k\times S^{m-k}$.  Let $f: S^m \mero S^k\times S^{m-k}$
be the associated meromorphic map, and take the map
\begin{equation}  \label{idf}
\id \times f: S^m\times S^m \mero S^m\times S^k\times S^{m-k}.
\end{equation}
By Proposition \ref{pull}, it suffices to
prove that the manifold $X=S^m\times S^k\times S^{m-k}$
is systolically free.  Let $\{g_j\}$ be a sequence of metrics on
$S^m\times S^k$ as in \eqref{systseq}. Let $h_j$ be a metric
on $S^{m-k}$ of volume $\vol_{m-k}(h_j)= \frac{\sys_m(g_j)}{\sys_k(g_j)}$.
Consider the metric $g_j\oplus h_j$ on $X$, and let $z$ be a
cycle representing a non-zero multiple of
$[S^k\times S^{m-k}]$.  Let $p: X\rightarrow S^{m-k}$ be the
projection to the last factor.  By the coarea inequality, we
obtain the following lower bound for the volume of
$z$ in $(X,g_j\oplus h_j)$:
\begin{equation}  \label{volm}
\vol_m(z)\geq \int_{(S^{m-k},h_j)} \vol_k(z\cap p^{-1}(x))\, dx\geq
\vol_{m-k}(h_j)\sys_k(g_j)=\sys_m(g_j),
\end{equation}
where the middle inequality uses intersection numbers for cycles and
transversality arguments in the context of maps of manifolds into $X$
(\cf Definition~\ref{def:systole} and \cite{BaK}, Lemma 6.1.).  Hence
\begin{equation}  \label{sysm}
\frac{\vol_{2m}(g_j\oplus h_j)}{\sys_m^2(g_j\oplus h_j)}=
\frac{\vol_{m-k} (h_j)\vol_{m+k}(g_j)}{\sys_m^2(g_j)}=
\frac{\vol_{m+k}(g_j)}{\sys_k(g_j)\sys_m (g_j)}
\xrightarrow[j\to\infty]{} 0,
\end{equation}
proving the systolic freedom of $X$ and therefore that of
$S^m\times S^m$.
\end{proof}

\section{Whitehead products and maps to wedges of spheres}
\label{sec:homotopy}

In this section, we establish some lemmas that will be needed in the
proof of Theorem~\ref{thm:mainCW}.  The main idea is to use
high-degree self-maps of the $m$-sphere $S^m$ as in R.~Thom~\cite{T}
to handle torsion in homotopy.  In what follows, we will denote a
wedge of $b$ copies of $S^m$ by $\vee^b S^m_r$, or simply $\vee S^m$.

A basic tool is the Hilton-Milnor theorem, which computes the homotopy
groups of a wedge of spheres in terms of the homotopy groups of
the factors.  In particular, it gives
the following splittings (see~\eg~\cite{Wh}):
\begin{subequations}  \label{eq:hilton}
\begin{alignat}{2}
\pi_k(\vee S^m_r)&=\bigoplus\nolimits_{r} \pi_k(S^m_r)
&\qquad & \text{ for }  k\le 2m-2 , \label{eq:hilton1}\\
\pi_{k}(\vee S^m_r)&=\bigoplus\nolimits_{r} \pi_{k}(S^m_r) \oplus
\bigoplus\nolimits_{r<s} \Z [e_r,e_s],
&\qquad & \text{ for }  k= 2m-1, \label{eq:hilton2}
\end{alignat}
\end{subequations}
where $[e_r,e_s]$ is the Whitehead product of the fundamental classes
of the corresponding spheres.

Another tool that we need is a formula of B.~Eckmann~\cite{E}
and G.~Whitehead~\cite{Wh1} on
the ``distributive law" in maps from spheres to spheres.
Let $\phi_q:S^m\to S^m$ be a map of degree $q$.
Let ${\phi_q}_{\sharp}: \pi_{k} (S^m) \to \pi_{k}(S^m)$ be the
induced map.  Then, if $x\in \pi_k(S^m)$,
\begin{subequations}  \label{eq:gwh}
\begin{alignat}{2}
{\phi_q}_{\sharp}(x)&= qx
&\qquad & \text{ for }  k\le 2m-2 , \label{eq:gwh1}\\
{\phi_q}_{\sharp}(x)&= qx + \binom{q}{2} H(x)\, [e,e]
&\qquad & \text{ for }  k= 2m-1,       \label{eq:gwh2}
\end{alignat}
\end{subequations}
where $H(x)\in \Z$ is the Hopf invariant of $x$, and
$e \in \pi_{m}(S^m)$ is the fundamental class, see \cite{Wh},~p.~537.

Let $K$ be a finite CW-complex of dimension $2m-1$.
Let $K^i$ be the $i$-skeleton of $K$.

\begin{lem} \label{lem1}
There exists a map $f:K \to \vee^b S^m$, where $b=\rank H_m(K/K^{m-1})$,
such that
\begin{equation}  \label{eq:inj}
f_*: H_m(K)\slash \tors (H_m(K)) \to H_m(\vee^b S^m)
{\text { is injective}}.
\end{equation}
\end{lem}

\begin{proof}   From the exact sequence of the pair $(K,K^{m-1})$
we see that the homomorphism $H_m(K)\rightarrow H_m(K/K^{m-1})$ is injective.
This reduces the problem to the case when $K$ is $(m-1)$-connected.
If $m=1$, there is a homotopy  equivalence $f:K\to \vee^b S^1$,
and we are done.  Thus, we may also assume that $m\ge 2$.

Let $H_m(K)\rightarrow \Z^{b}$, $b=b_m(K)$, be the quotient by the
torsion subgroup.  Since $K$ is $(m-1)$-connected, this homomorphism
is realized by a map to the corresponding Eilenberg-MacLane space,
$f:K\rightarrow K(\Z^{b},m)$.  Now $K(\Z^{b},m)=(K(\Z,m))^{b}$, and
$K(\Z,m)$ is obtained from $S^m$ by adding cells in dimensions $m+2$
and higher.   Hence, by the cellular approximation theorem, we get a map
$f_{m+1}:K^{m+1} \to \vee^{b} S^m$ by restricting to $(m+1)$-skeleta.
Clearly, this map satisfies \eqref{eq:inj}.
Thus, if $m=2$, we are done.

If $m\ge 3$, we extend the map $f_{m+1}$ to
$f=f_{2m}:K^{2m} \to \vee^{b} S^m$ by induction on skeleta,
proceeding as in Serre~\cite{S}, pp.~278 and 287--288.
Given $i$ with $m+1\leq i\leq 2m-2$, let
$f_{i}:K^i\to\vee^b S^m$ be a map satisfying \eqref{eq:inj}.
By \cite{S}, the group $\pi_i(S^m)$ is finite; let $q$ be its order.
Then, by  \eqref{eq:hilton1} and \eqref{eq:gwh1}, the map
$(\vee\phi_q)\circ f_{i}$ extends to a map $f_{i+1}:K^{i+1}\to\vee^b S^m$.
Clearly, $f_{i+1}$ satisfies \eqref{eq:inj}. This completes
the inductive step, and the proof.
\end{proof}

\begin{lem}\label{lem2}
A wedge of $m$-spheres, $m\ge 2$, admits a self-map
$\phi:\vee^b S^m\to \vee^b S^m$ such that
\begin{enumerate}
\item The map ${\phi}_{\sharp}: \pi_{2m-1}(\vee^b S^m)
\to \pi_{2m-1}(\vee^b S^m)$
has image contained in the subgroup generated by all
the Whitehead products; \label{item:pi}
\item The map $\phi_*: H_m(\vee^b S^m) \to H_m(\vee^b S^m)$
is injective.  \label{item:h}
\end{enumerate}
\end{lem}

\begin{proof}  For $b=1$, every Whitehead product is proportional to
$[e,e]$, where $e$ is the fundamental class of $S^m$.
Let $x\in \pi_{2m-1}(S^{m})$.  Then we can write
\begin{equation} \label{eq:hopf}
2x = s + H(x) [e,e]
\end{equation}
for some $s$ of finite order. Consider the map
$\phi=\phi_q: S^m \to S^m$, with $q$ even.
Substituting \eqref{eq:hopf} into \eqref{eq:gwh2}, we obtain:
\begin{align} \label{eq:whtoda}
{\phi}_{\sharp}(x) &=  \frac{q}{2} \left( s + H(x) [e,e] \right)
+ \binom{q}{2} H(x) [e,e] \\
&= \frac{q}{2} s + \frac{q^2}{2} H(x) [e,e].  \notag
\end{align}
Thus, it suffices to take
\begin{equation} \label{eq:q}
q = 2\left| \tors\left( \pi_{2m-1}(S^{m}) \right) \right|.
\end{equation}
Alternatively, a theorem of G. Whitehead \cite{Wh1} insures that
the subgroup of $\pi_{2m-1}(S^m)$ generated by Whitehead products
is precisely the kernel of the suspension homomorphism
$E:\pi_{2m-1}(S^m)\rightarrow\pi_{2m}(S^{m+1})$. Thus,
it suffices to pick $q=\left|\pi_{2m}(S^{m+1})\right|$ ,
which is less than or equal to the value from \eqref{eq:q}.

For $b>1$, we pick $\phi=\vee^b \phi_q$, with $q$ as in \eqref{eq:q}.
The splitting from \eqref{eq:hilton2} and an argument as above insure
that $\phi$ satisfies (\ref{item:pi}) and (\ref{item:h}).
\end{proof}

\begin{rem} \label{modd}
For $m$ odd, $m\ge 3$, we can actually choose $q$ so that
${\phi}_{\sharp}=0$, since, in that case, $\pi_{2m-1}(S^{m})$ is a
finite group, and all its elements have Hopf invariant $0$.
\end{rem}

\begin{rem} \label{projective}
Let $\F$ be $\C$, $\H$ or $\Ca$, and let $K=\F\P^2$ be the
corresponding projective plane.  With the usual decomposition into
$3$ cells for $K$, we have $K^{2m-1}=\F\P^1=S^m$, where $m=\dim(\F)$.
The smallest positive integer $q$ for which $\phi_q:S^m\to S^m$
satisfies conditions (\ref{item:pi}) and (\ref{item:h}) from
Lemma~\ref{lem2} can be computed explicitly in these examples.
Recall that $\pi_{2m-1}(S^m)=\Z\oplus T_m$, where $T_m$ is a finite
cyclic group, of order equal to $1$, $12$, or $120$ respectively
when $\F$ is $\C$, $\H$, or $\Ca$.  Let $a$ be the infinite order
generator defined by the Hopf map, and let $s$ be a generator of the
torsion part (taken to be $0$ when $\F=\C$).  A result of H.~Toda
(see \cite{Hu}) states that $[e,e]=2a\mp s$.  From this formula and
\eqref{eq:whtoda} we obtain for $q$ even:
\begin{equation} \label{eq:toda}
\phi_{q\sharp}(a) = \frac{q}{2}\, 2a+\binom{q}{2}\, [e,e]
           = \pm \frac{q}{2}\, s + \frac{q^2}{2}\, [e,e].
\end{equation}
Thus, the necessary and sufficient condition for
(\ref{item:pi}) and (\ref{item:h}) to hold is that
$q$ be a non-zero multiple of $2 \left| T_m\right|$.
\end{rem}

\section{Meromorphic maps to $S^m \times S^m$}
\label{sec:b1}

In this section, we prove Theorem~\ref{thm:mainCW} in the
particular case where $H_m(K/K^{m-1})=\Z$, by constructing
a meromorphic map $h:K \mero S^m\times S^m$.
The essential ingredients of the general case are already
present here, but the proof is more transparent in this
simpler situation.

\begin{proof}[Proof of Theorem~\ref{thm:mainCW} (particular case)]
Let $K$ be a finite, regular CW-complex of dimension $2m\ge 6$.
Assume $H_m(K/K^{m-1})=\Z$.

By Lemma~\ref{lem1}, there is a map $f:K^{2m-1}\to S^m$
such that $f_*: H_m(K^{2m-1}/K^{m-1}) \to H_m(S^m)$
is injective.  By Lemma~\ref{lem2},
there is a map $\phi=\phi^m_q: S^m\to S^m$ of degree $q\ne 0$ that maps
$\pi_{2m-1}(S^m)$ to the subgroup generated by the Whitehead
product $[e,e]$.  The map $\phi\circ f:K^{2m-1}\to S^m$ also maps
$\pi_{2m-1}(K^{2m-1})$ to this subgroup, while inducing
a monomorphism on $H_m$.

Now let $a_1$ and $a_2$ be the generators of $\pi_m(S^m \times S^m)$,
corresponding to the inclusions of the factors.
Recall that $S^m \times S^m = S^m \vee S^m \cup_{[a_1,a_2]} B^{2m}$.
Attaching an $(m+1)$-cell along the diagonal map
$(1,1):  S^m \to S^m \times S^m$,
we obtain the $2m$-dimensional regular CW-complex
\begin{equation}  \label{eq:wcw}
W = S^m \times S^m \cup_{a_1+a_2} B^{m+1}.
\end{equation}
Since $m\ge 2$, the complex $W$ satisfies condition~(\ref{mero1}) in
Definition~\ref{def:mero}.   Let $\a:S^m\rightarrow W$ be the composite
$S^m\xrightarrow{(1,0)} S^m \times S^m\hookrightarrow W$.
Then
\begin{alignat}{2} \label{eq:whiteprod}
\a([e,e]) &= [a_1,a_1]   &\qquad & \text{ since } \a(e)=a_1\\
         &= [a_1,-a_2]  &\qquad & \text{ since } a_1 + a_2 = 0
                           \text{ in } \pi_m(W)  \notag \\
         &= 0           &\qquad & \text{ since } [a_1,a_2]=0
                           \text{ in }\pi_{2m-1}(S^m\times S^m). \notag
\end{alignat}

Now let $h: K^{2m-1}\to W$ be the composite
$K^{2m-1}\xrightarrow{f} S^m\xrightarrow{\phi} S^m \xrightarrow{\a} W$.
By the above, $h_{\sharp}: \pi_{2m-1}(K^{2m-1}) \to \pi_{2m-1}(W)$
is the $0$ map.  Thus, $h$ extends over the $2m$-cells of $K$,
to a map $h: K\to W$.  Since $h$ is clearly injective on $H_m$,
we have defined a meromorphic map from $K$ to $S^m \times S^m$.
By Proposition~\ref{prop:snsn}, the manifold $S^m \times S^m$
is systolically free. Hence, by Proposition~\ref{pull}, the complex
$K$ is also systolically free.
\end{proof}

\section{Meromorphic maps to skeleta of products of spheres}
\label{sec:bgt1}

Before proving the general case of Theorem~\ref{thm:mainCW}, we
establish the systolic freedom of a model space by a ``long cylinder''
argument.

Let $X$ be a triangulated manifold of dimension $n$.  Let $A$ be the
$n$-skeleton of $X\times I$ where $I$ is an interval.  Then
$A=X\times\partial I\cup X^{n-1}\times I$.  Let $g_{+}$ and $g_{-}$ be
two metrics on $X$, and $g_{0}$ another metric dominating both $g_{+}$
and $g_{-}$.  Let $g$ be the metric on $A$ obtained by restricting the
metric $g_t\oplus dt^2$ of $X\times I$, where $I = [-L,L]$, with
$L=\ell+1>1$, and
\begin{equation} \label{eq:metric}
g_t =
\begin{cases}  g_{0}
&\qquad  \text{ if }  \left| t \right|\le \ell ,\\
(1-\lambda) g_{0} + \lambda g_{\pm}
&\qquad  \text{ if }  t =\pm ( \ell + \lambda ),  \text{ with }
0\le \lambda \le 1.
\end{cases}
\end{equation}

\begin{lem} \label{sublem}
For $k\geq 2$ and $\ell$ sufficiently large, we have
$\sys_k(g)\ge \beta$, where
\begin{equation} \label{eq:foo}
\beta= \frac{1}{2} \min\left( \sys_k(g_{+}), \sys_k(g_{-})\right).
\end{equation}
\end{lem}

\begin{proof}
Suppose $z$ is a non-bounding $k$-cycle in $A$ such that
$\vol_k(z)<\beta$.  Let $p: A\to [0,L]$ be the restriction of
the map $X\times I\to [0,L]$ given by
$(x,t)\mapsto \left| t \right|$.
The coarea inequality yields a point $t_0\in [0,\ell]$
such that the $(k-1)$-cycle $\gamma=z\cap p^{-1}(t_0)$ satisfies
\begin{equation}  \label{eq:vol}
\vol_{k-1}(\gamma)\le \frac{1}{\ell}\vol_{k}(z).
\end{equation}

By the isoperimetric inequality for cycles of small volume (\cite{G1},
Sublemma 3.4.B$'$) applied to $g_0|_{X^{n-1}}$, there is a constant
$C=C(g_0|_{X^{n-1}})$ with the following property: Every $(k-1)$-cycle
$\gamma$ in $X^{n-1}$ with $\vol(\gamma)<\frac{1}{C}$ bounds a
$k$-chain $D$ in $X^{n-1}$, of volume
\begin{equation}  \label{eq:volcyc}
\vol_k(D)\le C\vol_{k-1}(\gamma)^{\frac{k}{k-1}}.
\end{equation}
By choosing $\ell>\beta C$ we insure that the isoperimetric
inequality applies to $\gamma$.  Moreover, we need to choose $\ell$ so
that $\vol_k(D)<\beta$.  Thus we also require
$C\left(\frac{\beta}{\ell}\right)^{\frac{k}{k-1}}<\beta$, that is,
$\ell>\sqrt[\leftroot{1}\uproot{1} k]{\beta C^{k-1}}$.

Write $D=D_{-}+D_{+}$ where $D_{\pm}\subset X^{n-1}\times\{\pm t_0\}$.
Consider the decomposition of $z$ into a sum of cycles, $z=z_{-}+
z_{0} + z_{+}$, where $z_{-}=p^{-1}([0,t_0]) \cap (X\times
[-L,0])+D_{-}$, $z_{0}=p^{-1}([0,t_0])+D_{+}-D_{-}$, and
$z_{+}=p^{-1}([0,t_0]) \cap (X\times [0,L])-D_{+}$.  Now let
$\epsilon=0$, $+$, or $-$.  We have
\begin{equation}  \label{eq:sumcyc}
\vol_k(z_\epsilon) \le
\vol_k(z)+\vol_k(D)<\beta+\beta =
\min\left(\sys_k(g_+),\sys_k(g_-)\right).
\end{equation}
Hence $z_\epsilon$ is a boundary for every $\epsilon$ and so $[z]=0$.
The contradiction proves the lemma.
\end{proof}

\begin{lem} \label{lem4}
Let $B=(\times^c S^m)^{2m}$ be the $2m$-dimensional skeleton of
a product of $c$ copies of the $m$-sphere, $c\ge 2$, $m\ge 3$.
Then $B$ is systolically free.
\end{lem}

\begin{proof}
We choose the following representative $B_0$ in the homotopy class
of $B$.  Take the Cartesian product of the wedge $\vee^c S^m$ with
the wedge of $\binom{c}{2}$ intervals $I_{rs}=[0,L]$ for
sufficiently large $L=\ell+1$.  At the end of each interval, attach
a $2m$-cell along the Whitehead product, $[e_r,e_s]$, of the
fundamental classes of the spheres $S^m_r\times\{L\}$ and
$S^m_s\times \{L\}$ in $\vee^c S^m\times I_{rs}$:
\begin{equation}  \label{eq:model}
B_0=\left( \bigvee\nolimits^c S^m\right) \times
\left( \bigvee\nolimits_{r<s} I_{rs}\right) \cup
\bigcup\nolimits_{r<s} D_{rs}^{2m}.
\end{equation}

We precompose the projection $p:\vee I_{rs} \to [0,L]$ with the
projection to the second factor of $\vee S^m\times \vee I_{rs}$ and
extend it to a map $p: B_0\to [0,L]$ by setting $p(D_{rs}^{2m}) = L$.

Now we apply the argument of Lemma~\ref{sublem} with $n=2m$ to $B_0$,
as follows.  Identify each cell closure
$X_{rs}=\overline{D}^{\,2m}_{rs}$
with $S^m \times S^m$ by means of a diffeomorphism
$X_{rs}\xrightarrow{\rho_{rs}} S^m \times S^m$.  Then pull back a
sequence of free metrics $\{g_j\}$ on $S^m \times S^m$ provided by
Proposition \ref{prop:snsn} to obtain a sequence of free metrics
$\{\rho_{rs}^{*}(g_j)\}$ on $X_{rs}$, which play the role of the
metrics $g_\pm$ from Lemma \ref{sublem}.

Let $(X^{n-1},g_0)$ be the wedge of round spheres $S^m\vee S^m$ of
sufficiently big radius so that $g_0$ dominates all of the metrics
$\rho_{rs}^{*} (g_j)|_{X^{n-1}}$.  We obtain a lower bound for the
$m$-volume of a non-bounding $m$-cycle $z$ in $B_0$ by means of a
decomposition $z=z_0 + \sum_{r<s} z_{rs}$, where $z_{rs}$ is a cycle
in $\vee^c S^m \times I_{rs} \cup_{[e_r,e_s]} D^{2m}_{rs} \subset
B_0$.  This decomposition is obtained by the coarea inequality applied
to the projection $p: B_0\to[0,L]$.  This proves the systolic freedom
of $B_0$, and hence that of $B$, by Remark~\ref{rem:regCW}.
\end{proof}

\begin{proof}[Proof of Theorem~\ref{thm:mainCW} (general case)]
Let $K$ be a finite,
regular CW-complex of dimension $2m\ge 6$, with $\tors(H_m(K))=0$.
Let $f:K^{2m-1}\to \vee^b S^m$ be a map as in  Lemma~\ref{lem1},
and let $\phi: \vee^b S^m_r\to \vee^b S^m$ be a map
as in  Lemma~\ref{lem2}.  Then the composite
$\phi\circ f:K^{2m-1}\to \vee^b S^m$ maps
$\pi_{2m-1}(K^{2m-1})$ to the subgroup of $\pi_{2m-1}(\vee^b S^m)$
generated by the Whitehead products in the wedge, while inducing
a monomorphism on $H_m$.

Let $B$ be the $2m$-skeleton of a product of $c=2b$ copies of $S^m$,
as in Lemma \ref{lem4}.  We construct a meromorphic map from $\vee^{b}
S^m$ to $B$ much as in the case $b=1$.  Namely, let $\a_r: S^m_r\to
S^m_r\times S^m_{b+r}$ be the inclusion into the first factor.  Let
$W$ be the CW-complex obtained by attaching $(m+1)$-cells to $B$ along
the ``diagonals" $a_r+a_{b+r}$.  Define $\a: \vee^b S^m\to W$ to be
the composite
\begin{equation}  \label{alpha}
\a: \vee^b S^m_r\xrightarrow{\vee^b \a_i}\vee^b (S^m_r\times S^m_{b+r})
\hookrightarrow B \hookrightarrow W,
\end{equation}
and let $h=\a\circ \phi \circ f: K^{2m-1}\to W$.  Then $h_{\sharp}$ sends
$\pi_{2m-1}(K^{2m-1})$ to $0$.  Therefore, $h$ extends to a map $h: K\to W$.
Clearly, $h_*:H_m(K)\to H_m(W)$ is injective.  Thus,
we have defined a meromorphic map from $K$ to $X$.
By Lemma~\ref{lem4} and Proposition~\ref{pull}, $K$ is systolically free.
\end{proof}

\bigskip

\bibliographystyle{amsalpha}

\begin{thebibliography}{99}

\bibitem{Bk} I.~Babenko,
{\em Asymptotic invariants of smooth manifolds},
Russian Acad. Sci. Izv. Math. {\bf 41} (1993), 1--38.

\bibitem{BaK} I.~Babenko and M.~Katz,
{\em Systolic freedom of orientable manifolds},
Ann. Sci. \'{E}cole Norm. Sup. (1998, to appear); available at
{\ttfamily http://front.math.ucdavis.edu/math.DG/9707102}.

\bibitem{BeK} L.~B\'erard Bergery and M.~Katz,
{\em Intersystolic inequalities in dimension $3$},
Geom. Funct. Anal. {\bf 4} (1994), 621--632.

\bibitem{Ber1} M.~Berger,
{\em Du c\^ot\'e de chez Pu},
Ann. Sci. \'{E}cole Norm. Sup. {\bf 5} (1972), 1--44.

\bibitem{Ber2} \bysame,
{\em Systoles et applications selon Gromov},
S\'eminaire N.~Bourbaki, expos\'e 771, Ast\'erisque
{\bf 216} (1993), 279--310.

\bibitem{E} B.~Eckmann,
{\em Ueber die Homotopiegruppen von Gruppenra\"umen},
Comment. Math. Helv. {\bf 14} (1941), 234--256.

\bibitem{Gluck} H.~Gluck, D.~Mackenzie, and F.~Morgan, 
{\em Volume-minimizing cycles in Grassmann manifolds}, Duke
Math. J. {\bf 79} (1995), 335--404.

\bibitem{G1} M.~Gromov,
{\em Filling Riemannian manifolds}, J. Differential Geom.
{\bf 18} (1983), 1--147.

\bibitem{G2} \bysame,
{\em Systoles and intersystolic inequalities},
Actes de la table ronde de g\'eom\'etrie
diff\'erentielle en l'honneur de Marcel Berger (A.~Besse, ed.),
S\'eminaires et Congr\`es 1, Soci\'et\'e Math\'ematique de France,
1996, pp.~291--362.

\bibitem{G3} \bysame, 
{\em Metric structures for Riemannian and non-Riemannian spaces},
Progr. Math., vol. 152, Birkh\"auser, Boston, MA, 1998.

\bibitem{Hu} S.-T.~Hu,
{\em  Homotopy theory}, Pure Appl. Math.,
vol.~8, Academic Press, Boston, MA, 1959.

\bibitem{K1} M.~Katz,
{\em Counterexamples to isosystolic inequalities},
Geom. Dedicata {\bf 57} (1995), 195--206.

\bibitem{K2} \bysame,
{\em Systolically free manifolds}, Appendix~D in \cite{G3}.

\bibitem{LW} A.~Lundell and S.~Weingram,
{\em  The topology of CW complexes}, University Series in Higher Math.,
Van Nostrand Reinhold, New York, 1969.

\bibitem{P} C.~Pittet,
{\em Systoles on $S^1\times S^n$},
Differential Geom.\ Appl. {\bf 7} (1997), 139--142.

\bibitem{Sa}  P.~Sarnak,
{\em Extremal geometries}, Extremal Riemann surfaces
(J.~R.~Quine and P.~Sarnak, eds.), Contemp. Math.,
vol.~201, Amer. Math. Soc., Providence, RI, 1997, pp.~1--7.

\bibitem{S} J.-P.~Serre,
{\em  Groupes d'homotopie et classes de groupes abeliens},
Ann. of Math. {\bf 58} (1953), 258--294.

\bibitem{T} R.~Thom,
{\em Quelques propri\'et\'es globales des
vari\'et\'es diff\'erentiables},
Comment.\ Math.\ Helv. {\bf 28} (1954), 17--86.

\bibitem{Wh1} G.~W.~Whitehead,
{\em A generalization of the Hopf invariant},
Ann. of Math. {\bf 51} (1950), 192--237.

\bibitem{Wh} \bysame,
{\em  Elements of homotopy theory}, Grad. Texts
in Math., vol.~61, Springer-Verlag, New York, 1978.

\end{thebibliography}

\end{document}